\documentclass{article}

\usepackage{amsmath}
\usepackage{amsfonts}
\usepackage{amssymb}
\usepackage{amsthm}

\newtheorem{lem}{Lemma}[section]

\newtheorem{prop}[lem]{Proposition}
\newtheorem{thm}[lem]{Theorem}
\newtheorem{conj}[lem]{Conjecture}
\newtheorem{Defn}[lem]{Definition}
\newtheorem{Ex}[lem]{Example}

\newenvironment{defn}{\begin{Defn}\rm}{\end{Defn}}
\newenvironment{ex}{\begin{Ex}\rm}{\end{Ex}}

\newcommand{\ideal}[1]{\mathfrak{#1}}
\newcommand{\aaa}{\ideal{a}}

\newcommand{\m}{\ideal{m}}
\newcommand{\n}{\ideal{n}}
\newcommand{\p}{\ideal{p}}
\newcommand{\q}{\ideal{q}}
\newcommand{\rr}{\ideal{r}}
\newcommand{\s}{\ideal{s}}

\newcommand{\ol}{\overline}

\newcommand{\ti}{\tilde}

\newcommand{\x}{\mathbf{x}}
\newcommand{\y}{\mathbf{y}}

\newcommand{\rank}{\mathrm{rank}}

\newcommand{\Ht}{\mathrm{ht}\,}	

\newcommand{\spec}{\text{Spec}}	
	
\newcommand{\supp}{\text{Supp}}	
	
\newcommand{\len}{\mathrm{len}}
	
\newcommand{\ann}{\text{Ann}}

\begin{document}

\bibliographystyle{amsplain}
\begin{center}
\large \textbf{Multiplicities and a Dimension Inequality for Unmixed Modules}

\medskip
\normalsize Sean Sather-Wagstaff\footnote{This research was partially
conducted by the author for the Clay Mathematics Institute.}

\medskip
Department of Mathematics, University of Illinois, 273 Altgeld Hall,\\
1409 West Green Street, Urbana IL, 61801\\
Email: {\tt ssather@math.uiuc.edu}
\end{center}

\bigskip
\small 
We prove the following result, which is motivated by
the recent work of Kurano and Roberts on Serre's positivity 
conjecture.
Assume that $(R,\m)$ is a local ring with 
finitely-generated module $M$ such that $R/\ann(M)$ is quasi-unmixed 
and contains a field,
and that $\p$ and $\q$ are prime ideals in the support of $M$ such 
that $\p$ is analytically unramified, $\sqrt{\p+\q}=\m$ and
$e(M_{\p})=e(M)$.  Then
\[ \dim(R/\p)+\dim(R/\q)\leq \dim(M). \]
We also prove a similar theorem in a special case of mixed 
characteristic.
Finally, we provide several examples to explain our assumptions as 
well as an 
example of a noncatenary, 
local domain $R$ with prime ideal $\p$ such that $e(R_{p})>e(R)=1$.

\normalsize

\section{Introduction}
Throughout this paper, all rings are commutative and Noetherian, and 
all modules are finitely 
generated and unital.

Let $(R,\m)$ be a local Noetherian ring of dimension $d$, and let $M$ 
and $N$ be finitely 
generated $R$-modules such that $M$ has finite projective dimension 
and $M\otimes_R N$ 
is a module of finite length.  Serre~\cite{serre:alm} defined the 
\textit{intersection multiplicity} 
of $M$ and $N$ to be 
\[ \chi (M,N)=\sum_{i=0}^d (-1)^i \mbox{length} 
(\mbox{Tor}^R_i(M,N))\]  
and considered the following properties when $R$ is regular.
\begin{enumerate}
\item $\dim (M)+\dim (N)\leq \dim (R)$.
\item (Nonnegativity) $\chi (M,N) \geq 0$.
\item (Vanishing) If $\dim (M)+\dim (N)< \dim (R)$, then $\chi (M,N) 
= 0$.
\item (Positivity) If $\dim (M)+\dim (N)=\dim (R)$, then $\chi 
(M,N)>0$.
\end{enumerate}
Serre was able to verify the first statement for any regular local 
ring and the others 
in the case when $R$ is unramified.  Since 
$\chi (M,N)$ has many of the characteristics we desire from an 
intersection multiplicity 
(for example, B\'ezout's Theorem holds), it was reasonable to suppose 
that these further 
properties are satisfied over an arbitrary regular local ring.   The 
results were left unproved 
for ramified rings.  

The vanishing conjecture was proved about ten years ago by Gillet and 
Soul\'{e}~\cite{gillet:knmi}
and independently by Roberts~\cite{roberts:vimpc} using 
$K$-theoretic methods. Gabber proved 
the nonnegativity conjecture 
recently~\cite{berthelot:ava,hochster:nimrrlr,roberts:rdsmcgpnc} 
using a theorem of de Jong~\cite{dejong:ssa}.  The positivity 
conjecture remains open in the ramified
case.  Kurano and Roberts have 
proved the following using methods introduced by Gabber.

\begin{thm} \label{thm:0}
\textnormal{(\cite{kurano:pimsppi} Theorem 3.2)}
Assume that $(R, \ideal{m})$ is a regular local ring which either 
contains a field or is ramified.  
Also, assume that $\ideal{p}$ and $\ideal{q}$ are prime ideals in $R$ 
such that 
$\sqrt{\ideal{p}+\ideal{q}}=\ideal{m}$ and 
$\dim(R/\ideal{p})+\dim(R/\ideal{q})=\dim R$.  If 
$\chi(R/\p,R/\q)>0$ then 
\begin{equation}
\ideal{p}^{(n)}\cap \ideal{q}\subseteq \ideal{m}^{n+1} \qquad\text{ 
for all } n> 0 \label{eqn:0}
\end{equation}
where $\p^{(n)}$ denotes the $n$th symbolic power of $\p$.
\end{thm}

As a result, they conjectured that (\ref{eqn:0}) should hold for all 
regular local rings.  

\begin{conj} \label{conj:1}
Assume that $(R, \ideal{m})$ is a regular local ring 
and that $\ideal{p}$ and $\ideal{q}$ are prime ideals in $R$ such 
that 
$\sqrt{\ideal{p}+\ideal{q}}=\ideal{m}$ and 
$\dim(R/\ideal{p})+\dim(R/\ideal{q})=\dim R$.  Then 
$\ideal{p}^{(n)}\cap \ideal{q}\subseteq \ideal{m}^{n+1}$ for all $n> 
0$.
\end{conj}

We study Conjecture~\ref{conj:1}, as a verification of this 
conjecture could introduce new tools
to apply to the positivity conjecture.

For any local ring $(A,\n)$ with finite module $M$ let $e(M)$ denote 
the Samuel multiplicity of $M$ 
with respect to the ideal $\n$.  
(For the definition of the Hilbert-Samuel multiplicity, see 
Definition~\ref{def:mult} below.)
It is straightforward to 
verify that, if $R$ is a regular local ring with prime ideal $\p$ and 
$0\neq f\in\p$, then 
$e(R_{\p}/(f))=m$ if and only if 
$f\in\p^{(m)}\smallsetminus\p^{(m+1)}$.  Thus, 
Conjecture~\ref{conj:1} may be rephrased as the following.

\bigskip
\noindent\textbf{Conjecture~\ref{conj:1}$'$}
\textit{Assume that $(R, \ideal{m})$ is a regular local ring 
and that $\ideal{p}$ and $\ideal{q}$ are prime ideals in $R$ such 
that $\sqrt{\ideal{p}+\ideal{q}}=\m$.  If there exists $0\neq 
f\in\p\cap\q$ such that $e(R_{\p}/(f))=e(R/(f))$, then 
$\dim(R/\p)+\dim(R/\q)\leq \dim(R)-1$.}

\bigskip

Conjecture~\ref{conj:1}$'$ motivates the following generalization.

\begin{conj} \label{conj:main}
Assume that $(R,\m)$ is a  local ring with finitely-generated module 
$M$ such 
that $R/\ann(M)$ is quasi-unmixed, 
and $\p$ and $\q$ are prime ideals in the support of $M$ such that 
$\p$ is analytically unramified, $\sqrt{\p+\q}=\m$ and 
$e(M_{\p})=e(M)$.  Then
\[ \dim(R/\p)+\dim(R/\q)\leq \dim(M). \]
\end{conj}

In Examples~\ref{ex:eq}--\ref{ex:morenagata} we give examples showing 
the necessity of the
assumptions of Conjecture~\ref{conj:main}.  

In a previous paper (\cite{sather:dicmr}) we considered the case 
$M=R$ and $R$ is Cohen-Macaulay.
In Conjecture~\ref{conj:main}, if $\dim(M)=\dim(R)$ and $R$ is 
quasi-unmixed, then the Associativity Formula for multiplicities
tells us that the condition $e(M_{\p})=e(M)$ is slightly weaker than 
the
condition $e(R_{\p})=e(R)$.   The condition $e(R_{\p})=e(R)$ depends
on all the minimal primes of $R$, while the condition $e(M_{\p})=e(M)$
depends only on the minimal primes of $R$ which are in the support of
$M$.
In Example~\ref{ex:more} we give an example demonstrating this.

The following is  the main result of this paper, in which we prove 
Conjecture~\ref{conj:main} in the case where $R/\ann(M)$ contains a 
field.  

\bigskip
\noindent\textbf{Theorem~\ref{thm:main}.}
\textit{Assume that $(R,\m)$ is a local ring
and that $M$ is an $R$-module such that $R/\ann(M)$ is 
quasi-unmixed and contains a
field.  Let $\p$ and $\q$ be prime ideals of $R$ in the support of 
$M$ such that $\p$ is analytically unramified,
$\sqrt{\p+\q}=\m$, and $e(M_{\p})=e(M)$.  Then 
$\dim(R/\p)+\dim(R/\q)\leq \dim(M)$.}  

\bigskip

Note that, if $R$ is assumed to be excellent, then the condition that 
$\p$ is analytically unramified is automatically satisfied, and 
$R/\ann(M)$ is quasi-unmixed if and only if it is equidimensional.

Here we give a sketch of the proof of the main theorem.  By passing 
to the quotient $R/\ann(M)$, we may assume that
$R$ is quasi-unmixed and contains a
field, and that $\supp(M)=\spec(R)$.   By using the Associativity 
Formula for multiplicities (see Section~\ref{sec:defs}), we may 
replace the module $M$ with the ring $R$.  By passing to the ring 
$R[X]_{\m R[X]}$,
we may assume that the residue field of $R$ is infinite.  The fact 
that $\p$ is analytically unramified allows us to pass to the
completion of $R$ so that we may assume that $R$ is complete and 
equidimensional and contains an infinite field $K$.
To prove the complete case, we construct an injection
$R'=K[\![X_1,\ldots,X_r]\!]\rightarrow R$
where $r=\dim(R)$, and
we let $\p'=\p\cap R'$ and $\q'=\q\cap R'$.  To prove the result, it 
suffices to show that
$\p'+\q'$ is primary to the maximal ideal of $R'$.  This is 
sufficient because $R'$ is regular and Serre's result shows 
that
\begin{align*}
\dim(R/\p)+\dim(R/\q) & = \dim(R'/\p')+\dim(R'/\q') \\
 & \leq\dim(R')=\dim(R).
\end{align*}
To prove that $\p'+\q'$ is primary to the maximal ideal of $R'$, it 
suffices to show that $\p$ is the unique prime
ideal of $R$ which contracts to $\p'$ in $R'$.  The desired
uniqueness follows from our key lemma, which is essentially
a corollary of a standard formula in multiplicity theory.

\bigskip
\noindent\textbf{Lemma~\ref{lem:1}.}
\textit{Assume that $R$ is an equidimensional ring containing a 
regular local ring $(R',\m')$ such that the extension 
$R'\rightarrow R$
is module finite.  Let $\p$ be a prime ideal of $R$ and let 
$\p'=\p\cap R'$.
Assume that $e(R_{\p})=\rank_{R'}(R)$.  Then $\p$ is the unique prime 
ideal of $R$ contracting
to $\p'$ in $R'$ and $\kappa(\p')\cong \kappa(\p)$.}

\bigskip

In addition, we use similar methods to prove the following theorem in 
mixed characteristic.

\bigskip
\noindent\textbf{Theorem~\ref{thm:mixed}}
\textit{Assume that $(R,\m)$ is a  local ring and that $M$ is an 
$R$-module such that $R/\ann(M)$ is quasi-unmixed of mixed 
characteristic $p$ and that $p$ is part of a system of parameters of 
$R/\ann(M)$ which generates a reduction of the maximal ideal of 
$R/\ann(M)$.  Let $\p$ and $\q$ be prime ideals of $R$ in the support 
of $M$ such that $\p$ is analytically unramified,
$\sqrt{\p+\q}=\m$, and $e(M_{\p})=e(M)$.  Then 
$\dim(R/\p)+\dim(R/\q)\leq \dim(M)$.}

\bigskip

As noted above, Serre proved Conjecture~\ref{conj:main} in the case 
where $R$ is regular and $M=R$.   
Many attempts have been made to generalize Serre's result to the
nonregular situation.  It can be shown easily that, if one drops the 
assumption of regularity on the ring, then one
must assume that the objects under investigation have additional 
properties which would be automatic if the ring
were regular.  In our conjectures and results, we assume that the 
module $M$ satisfies the property $e(M_{\p})=e(M)$.  
When the ring is regular and $M$ has positive rank, this is automatic 
because the localization $R_{\p}$ is also regular so that $e(R)=1$ and
$e(R_{\p})=1$, and if $r=\rank_R(M)$ then
$r=\rank_{R_{\p}}(M_{\p})$ and $e(M)=re(R)=r=re(R_{\p})=e(M_{\p})$.  
In a famous conjecture, Peskine and Szpiro focus on the finiteness of 
projective
dimensions over regular local rings.  

\begin{conj} \label{conj:3}
\textnormal{(\cite{peskine:dpfcl})}
Assume that $(R,\m)$ is a local ring 
and $\p$ and $\q$ are prime ideals in $R$ such that $\p$ has finite 
projective dimension and $\sqrt{\p+\q}=\m$.  Then 
$\dim(R/\p)+\dim(R/\q)\leq\dim(R)$.
\end{conj}

We include this here because it places Conjecture~\ref{conj:main} in 
a second context:  not only is 
our conjecture motivated by a consequence of positivity, but also it 
is a generalization of Serre's dimension 
inequality for regular local rings.
We note that, in the nonregular situation, 
Conjectures~\ref{conj:main} and~\ref{conj:3} are not 
comparable:  c.f., \cite{sather:dicmr}.

In Section~\ref{sec:defs} we present definitions and background 
results.  In Section~\ref{sec:main} we prove the key
lemma and the main theorem.  In Section~\ref{sec:mixed} we prove the 
mixed characteristic result mentioned above.  
In Section~\ref{sec:orig} we briefly consider Conjecture~\ref{conj:1} 
in light of Theorems~\ref{thm:main} and~\ref{thm:mixed}.  
In Section~\ref{sec:ex} we present several examples.

\section{Definitions and Background Results} \label{sec:defs}

The Hilbert-Samuel multiplicity shall play a central role in our 
work.  For the sake of clarity 
we specify which multiplicity 
we are considering.  

\begin{defn} \label{def:mult}
Assume that $(R,\m)$ is a local ring and that $M$ is an $R$-module of 
dimension $d$.
Let $\aaa$ be an ideal of $R$ such that
$\sqrt{\aaa}=\m$.   For $n\gg 0$ the 
Hilbert function
$H_{\aaa,M}(n)=\len(M/\aaa^{n+1}M)$
is a polynomial in $n$ 
of degree $d$ with rational coefficients.
If $e_d$ is the leading coefficient of this polynomial then the 
\textit{Hilbert-Samuel multiplicity} of 
$\aaa$ on $M$ is the positive integer $e_R(\aaa,M)=d! e_d$.  We will 
write $e(\aaa,M)$ instead of $e_R(\aaa,M)$
if doing so causes no
confusion.  We denote $e(\m,M)$ by $e(M)$.  
\end{defn}

Recall that the Hilbert-Samuel 
multiplicity satisfies the following \textit{Associativity Formula}
\[ e(M)=\sum_{\p}\len(M_{\p})e(R/\p) \]
where the sum is taken over all prime ideals $\p$ of $R$ such that 
$\dim(R/\p)=\dim(M)$.  Because we need only
take the sum over such prime ideals which are also in the support of 
$M$, this sum is finite.

Furthermore, recall the following formula (c.f., 
Nagata~\cite{nagata:lr} (23.1)).
Assume $(R',\m')$ is a local ring contained in a semilocal ring 
$(R,\m_1,\ldots,\m_n)$ such that each $\m_i\cap R'=\m'$.
Assume also that $R$ is a finite $R'$-module.  Then
\begin{equation} 
e_{R'}(\m',R)=\sum_i [R/\m_i:R'/\m'] e_{R_{\m_i}}(\m' 
R_{\m_i},R_{\m_i}) \label{eqn:mult}
\end{equation}
where the sum is taken over all indices $i$ such that 
$\Ht(\m_i)=\dim(R')$.

The following theorem tells us that, under certain circumstances, the 
Hilbert-Samuel multiplicity
is well-behaved under localizations.  We say that a local ring $R$ is 
\textit{analytically
unramified} if its completion has no nonzero nilpotents.  We say that 
a 
prime ideal $\p$ in a local ring $R$ is \textit{analytically
unramified} if the quotient $R/\p$ is analytically unramified.

\begin{thm} \label{thm:loc}
\textnormal{(\cite{nagata:lr} (40.1))}
Let $\p$ be a prime ideal of a local ring $R$.  If $\p$ is 
analytically
unramified and if $\Ht(\p)+\dim(R/\p)=\dim(R)$, then $e(R_{\p})\leq 
e(R)$.
\end{thm}

Note that, if $R$ is catenary and equidimensional then the condition 
$\Ht(\p)+\dim(R/\p)=\dim(R)$ is 
automatically satisfied.  (A ring $R$ is 
\textit{equidimensional} if, for every 
minimal prime $\rr$ of $R$,
$\dim(R/\rr)=\dim(R)$.)
Theorem~\ref{thm:loc} gives our motivation for the assumption ``$\p$ 
is analytically unramified'' 
in Conjecture~\ref{conj:main}.  If we do not assume that the 
Hilbert-Samuel multiplicity is 
well-behaved with respect to localization, then there is no reason to 
suspect that the assumption 
$e(M_{\p})=e(M)$ is meaningful.  
Example~\ref{ex:nagata1} is an example of a local domain with prime 
ideal $\p$ which is analytically 
unramified such that $\Ht(\p)+\dim(R/\p)\neq\dim(R)$ and 
$e(R_{\p})>e(R)$.  At this time, we do not 
know of an example of an equidimensional, catenary local ring $R$ 
with prime ideal $\p$ which is 
analytically ramified where either Conjecture~\ref{conj:main} or the 
conclusion of 
Theorem~\ref{thm:loc} fail to hold.  

The following result tells us that multiplicities behave well under 
certain flat extensions.

\begin{prop} \label{prop:flat}
\textnormal{(Herzog~\cite{herzog:omlr} Lemma 2.3)}
Assume that $R\rightarrow \ti{R}$ is a flat local homomorphism of 
local rings $(R,\m)$ and
$(\ti{R},\ti{\m})$ and that $\m\ti{R}=\ti{\m}$.  Then 
$e(\ti{R})=e(R)$.
\end{prop}

\section{The Main Theorem} \label{sec:main}

The key lemma for our main theorem is the following.  
For an
integral domain $A$ let $Q(A)$ denote the field of fractions of $A$.
For a prime ideal $\p$ of a ring $R$, let $\kappa(\p)=Q(R/\p)$.

\begin{lem} \label{lem:1}
Assume that $R$ is an equidimensional ring containing a 
regular local ring $(R',\m')$ such that the extension 
$R'\rightarrow R$
is module finite.  Let $\p$ be a prime ideal of $R$ and let 
$\p'=\p\cap R'$.
Assume that $e(R_{\p})=\rank_{R'}(R)$.  Then $\p$ is the unique prime 
ideal of $R$ contracting
to $\p'$ in $R'$ and $\kappa(\p')\cong \kappa(\p)$.  
\end{lem}

\begin{proof}
The fact that the extension $R'\rightarrow R$ is module finite and 
injective implies
that $R$ is semilocal and dominates $R'$.
Let $r=\rank_{R'}(R)$ so that 
$r=\rank_{R'_{\p'}}(R_{\p'})=\rank_{R'_{\p'}}(R\otimes_{R'}R'_{\p'})$.
By eq.\ (\ref{eqn:mult}) and Matsumura~\cite{matsumura:crt} Theorem 
14.8,
\[ r=r e(R'_{\p'})=e(\p' R'_{\p'},R_{\p'})=\sum_{\q} 
[\kappa(\q):\kappa(\p')] e(\p' R_{\q},R_{\q}) \]
where the sum is taken over all prime ideals $\q$ of $R$ which 
contract to $\p'$ and such that
$\Ht(\q)=\Ht(\p')$.  Because the extension
$R'\rightarrow R$ is finite, $R'$ is integrally closed and $R$ is 
equidimensional, the going-up
and going-down properties hold for the extension.  In particular, any 
prime ideal $\q$ of $R$ which
contracts to $\p'$ automatically satisfies the condition 
$\Ht(\q)=\Ht(\p')$.  Our assumption
$e(R_{\p})=r$ implies that
\[ r=\sum_{\q} [\kappa(\q):\kappa(\p')] e(\p' R_{\q},R_{\q})\geq e(\p'
R_{\p},R_{\p}) \geq e(R_{\p})=r \]
where the sum is taken over all prime ideals $\q$ such that $\q\cap 
R'=\p'$.
Therefore, we have equality at each step.  The only way this can be 
true is if $\p$ is the unique 
such prime and $[\kappa(\p):\kappa(\p')]=1$.  This is the desired 
conclusion.
\end{proof}

The following is our main theorem.  Recall that a local ring is 
\textit{quasi-unmixed} (or \textit{formally equidimensional}) if its 
completion is equidimensional.  This is equivalent to the ring being 
universally catenary and equidimensional by 
Ratliff~\cite{ratliff:qldafccpiII}.

\begin{thm} \label{thm:main}
Assume that $(R,\m)$ is a  local ring
and that $M$ is an $R$-module such that $R/\ann(M)$ is 
quasi-unmixed and contains a
field.  Let $\p$ and $\q$ be prime ideals of $R$ in the support of 
$M$ such that $\p$ is analytically unramified,
$\sqrt{\p+\q}=\m$, and $e(M_{\p})=e(M)$.  Then 
$\dim(R/\p)+\dim(R/\q)\leq \dim(M)$.
\end{thm}

As we noted above, if $R$ is assumed to be excellent, then the 
condition that $\p$ is analytically unramified is automatically 
satisfied, and $R/\ann(M)$ is 
quasi-unmixed if and only if it is equidimensional.

\begin{proof}
Step 1.  By passing to the quotient $R/\ann(M)$, we reduce to the 
case where $R$ is quasi-unmixed and contains a field and 
$\supp(M)=\spec(R)$.  
The fact that $\p$ and $\q$ are in the support of $M$ implies that 
$\ann(M)\subseteq\p\cap\q$.
Let $R'=R/\ann(M)$ with $\m'=\m R'$, $\p'=\p R'$ and $\q'=\q R'$.
Our assumptions imply that $R'$ is a quasi-unmixed local ring
which contains a field, and $\supp_{R'}(M)=\spec(R')$.  
The ideals $\p'$ and $\q'$ are primes in the support of $M$
such that $\p'$ is analytically unramified and $\sqrt{\p'+\q'}=\m'$.  
Also, $e_{R'}(M)=e_R(M)$ because the Hilbert polynomials of $M$
over $R$ and $R'$ are the same.  Similarly, 
$e_{R'_{\p'}}(M_{\p'})=e_{R_{\p}}(M_{\p})$ so that
$e_{R'_{\p'}}(M_{\p'})=e_{R'}(M)$. Thus, 
the case where $R$ is quasi-unmixed and contains a field and 
$\supp(M)=\spec(R)$ implies that
\[ \dim(R/\p)+\dim(R/\q)=\dim(R'/\p')+\dim(R'/\q')\leq \dim(M) \]
giving the desired result.

Step 2.  We reduce to the case where $M=R$.  Because we are assuming 
now that $\supp(M)=\spec(R)$
and that $R$ is quasi-unmixed, the Associativity Formula tells us that
\begin{equation}
\sum_{\rr}\len(M_{\rr})e(R/\rr)=e(M)=e(M_{\p})
        =\sum_{\rr\subseteq\p} \len(M_{\rr})e(R_{\p}/\rr_{\p}) 
\label{eqn:x}
\end{equation}
where the first sum is taken over all minimal primes of $R$ and the 
second sum is taken over all
minimal primes of $R$ which are contained in $\p$.  The fact that each
$R/\rr$ is 
quasi-unmixed and that $\p/\rr$ is analytically unramified for all 
such $\rr$ contained in $\p$
implies that $e(R_{\p}/\rr_{\p})\leq e(R/\rr)$ by 
Theorem~\ref{thm:loc}.  Thus,
eq.\ (\ref{eqn:x}) implies that every minimal prime $\rr$ of $R$ is 
contained in $\p$ and
$e(R_{\p}/\rr_{\p})= e(R/\rr)$.  Thus, 
\[ 
e(R)=\sum_{\rr}\len(R_{\rr})e(R/\rr)=\sum_{\rr}\len(R_{\rr})e(R_{\p}/\rr_{\p})=e(R_{\p}). 
\]
Of course, $\p$ and $\q$ are in the support of $R$, so that, if we 
know the result for $M=R$
then
\[ \dim(R/\p)+\dim(R/\q)\leq\dim(R)=\dim(M).\]

Step 3.  By passing to the ring $R(X)=R[X]_{\m R[X]}$, we reduce to
the case where $R$ has infinite residue field.  The ring 
$R(X)=R[X]_{\m R[X]}$ is quasi-unmixed with infinite residue
field $R/\m(X)$.  Because the extension $R\rightarrow R(X)$ is 
flat and local, and $\m$ extends to the maximal ideal of $R(X)$,
Proposition~\ref{prop:flat} implies that $e(R(X))=e(R)$.  For any
ideal $I$ of $R$ let $I(X)=I R(X)$.  Then $\p(X)$ and $\q(X)$ are 
prime ideals of $R(X)$ such
that $\sqrt{\p(X)+\q(X)}=\m(X)$.  By~\cite{nagata:lr} (36.8), $\p(X)$ 
is analytically unramified.
The extension $R_{\p}\rightarrow R(X)_{\p(X)}$ is faithfully
flat and $\p_{\p}$ extends to the maximal ideal of $R(X)_{\p(X)}$, so 
that 
$e(R(X)_{\p(X)})=e(R_{\p})$.    If we know the result for rings with 
infinite residue field, then
\begin{align*} 
\dim(R/\p)+\dim(R/\q)&=\dim(R(X)/\p(X))+\dim(R(X)/\q(X))\\
&\leq\dim(R(X))=\dim(R). 
\end{align*}

Step 4.  By passing to the completion $\hat{R}$ of $R$, we reduce to 
the case where $R$ is complete and equidimensional.  Let $\hat{\p}$ 
and $\hat{\q}$ be prime ideals of $\hat{R}$ which are minimal over 
$\p\hat{R}$ and
$\q\hat{R}$, respectively.  Then $\dim(\hat{R}/\hat{\p})=\dim(R/\p)$ 
and similarly for 
$\hat{R}/\hat{\q}$.  The ring $\hat{R}$ is quasi-unmixed.  Also,
\[ 
\hat{\m}\supseteq\sqrt{\hat{\p}+\hat{\q}}\supseteq\sqrt{\p\hat{R}+\q\hat{R}}=\m\hat{R}
        =\hat{\m} \]
so that $\sqrt{\hat{\p}+\hat{\q}}=\hat{\m}$.    The extension 
$R\rightarrow\hat{R}$ is flat and local
and the extension of $\m$ into $\hat{R}$ is the maximal ideal 
$\hat{\m}$, so that
$e(\hat{R})=e(R)$ by Proposition~\ref{prop:flat}.  The fact that $\p$ 
is analytically unramified implies that
$\p\hat{R}_{\hat{\p}}=\hat{\p}\hat{R}_{\hat{\p}}$, so
that $e(\hat{R}_{\hat{\p}})=e(R_{\p})$.  Thus, the complete, 
equidimensional case implies that
\[ 
\dim(R/\p)+\dim(R/\q)=\dim(\hat{R}/\hat{\p})+\dim(\hat{R}/\hat{\q})\leq\dim(\hat{R})
        =\dim(R). \]

Step 5.  We prove the case where $R$ is complete and equidimensional 
with infinite residue field, and $M=R$.
By Bruns and Herzog~\cite{bruns:cmr} Proposition 4.6.8, there exists 
a system of parameters
$\x=x_1,\ldots,x_r$ of $R$ which generates a reduction ideal of 
$\m$.  

Let $K$ denote the residue field of $R$, which we may assume is 
contained in $R$, as $R$
is complete and contains a field.  Let $R'$ denote the power series 
ring 
$K[\![ X_1,\ldots,X_r]\!]$.  The natural map $R'\rightarrow R$
given by $X_i\mapsto x_i$ is injective and makes $R$ into a finite 
$R'$-module.
By~\cite{matsumura:crt} Theorem 14.8 and~\cite{bruns:cmr} Lemma 4.6.5
\[ \rank_{R'}(R)=\rank_{R'}(R) 
e(R')=e(\m'R,R)=e((\x)R,R)=e(R)=e(R_{\p}).\]
Let $e=e(R)$, $\p'=\p\cap R'$ and $\q'=\q\cap R'$. 
The fact that the extension $R'/\p'\rightarrow R/\p$ is injective and
module-finite implies that $\dim(R'/\p')=\dim(R/\p)$, and similarly, 
$\dim(R'/\q')=\dim(R/\q)$.

To show that $\dim(R/\p)+\dim(R/\q)\leq\dim(R)$, it suffices to show 
that $\sqrt{\p'+\q'}=\m'$,
as Serre's intersection theorem for regular local rings implies that
\begin{align*}
\dim(R/\p)+\dim(R/\q) & = \dim(R'/\p')+\dim(R'/\q') \\
 & \leq\dim(R')=\dim(R) 
\end{align*}
as desired.

The fact that 
$e(R_{\p})=e(R)=\rank_{R'}(R)$ with the previous lemma implies that
$\p$ is the unique prime ideal of $R$ contracting to $\p'$ in $R'$.  
Let $\s'$ be a prime ideal of $R'$ containing $\p'+\q'$.  It suffices 
to
show that $\s'=\m'$.  By the going-up property, there is a prime 
ideal $\s$ of $R$ containing $\q$ such that $\s\cap R'=\s'$.  By the 
going-down property, there is a prime ideal $\p_1$ of $R$ contained 
in $\s$ such that $\p_1\cap R'=\p'$.  By our uniqueness statement, 
$\p_1=\p$ so that $\s$ contains $\p+\q$.  Therefore, $\s=\m$ and 
$\s'=\s\cap R'=\m\cap R'=\m'$, as desired.
This completes the proof.
\end{proof}

\section{A Theorem in Mixed Characteristic} \label{sec:mixed}

The following result in mixed characteristic is slightly more 
restrictive than Theorem~\ref{thm:main}, but is rather interesting.  
For brevity, we say that a system of parameters $x_1,\ldots,x_d$ of a 
local ring $(R,\m)$ is a \textit{reductive system of parameters} if 
$(\x)R$ is a reduction of $\m$.  If the residue field of $R$ is 
infinite and the sequence $y_1,\ldots,y_r$ generates a reduction of 
$\m$, then Northcott and Rees~\cite{northcott:rilr} Theorem 4.2 tells 
us that the ideal $(\y)R$ is a minimal reduction of $\m$ if and only 
if $r=\dim(R)$.  That is, a system of parameters is reductive if and 
only if it generates a minimal reduction of $\m$.

\begin{thm} \label{thm:mixed}
Assume that $(R,\m)$ is a local ring and that $M$ is an $R$-module 
such that $R/\ann(M)$ is quasi-unmixed of mixed characteristic $p$ 
and that $p$ is part of a reductive system of parameters of 
$R/\ann(M)$.  Let $\p$ and $\q$ be prime ideals of $R$ in the support 
of $M$ such that $\p$ is analytically unramified,
$\sqrt{\p+\q}=\m$, and $e(M_{\p})=e(M)$.  Then 
$\dim(R/\p)+\dim(R/\q)\leq \dim(M)$.
\end{thm}

Again, we note that if $R$ is excellent, then the assumption that 
$\p$ is analytically unramified is automatically satisfied and that 
$R/\ann(M)$ is quasi-unmixed if and only if it is equidimensional.

\begin{proof}
We follow the same steps as in the proof of Theorem~\ref{thm:main}.  
Steps 1--4 are independent of the assumption on the characteristic of 
the ring.  To verify that we may assume that $R$ is complete, 
equidimensional and mixed characteristic $p$ with infinite residue 
field and that $p$ is part of a reductive system of parameters of 
$R$, it suffices to show that this final property passes through each 
of the steps.  Steps 1 and 2 are trivial.  For Steps 3 and 4, we note 
that each extensions $R\rightarrow R(X)$ and $R\rightarrow \hat{R}$ 
are flat, local extensions such that the extension of the maximal 
ideal of $R$ into the extension ring is the maximal ideal of the 
extension ring.  It is straightforward to show that, in this 
situation, a reductive system of parameters $\x$ of $R$  extends to a 
reductive system of parameters of the extension ring.  Therefore we 
may pass to the completion.

Step 5.  Let $p=x_1,x_2,\ldots,x_d$ be a reductive system of 
parameters of $R$.  The fact that $p$ is part of a system of 
parameters for $R$ implies that $p$ is contained in no minimal prime 
of $R$.  In particular, $R$ has characteristic 0.
The fact that $R$ is complete then implies that $R$ has a coefficient 
ring $(V,pV)$ which is a complete discrete valuation ring contained 
in $R$.  Let $R'= V[\![X_2,\ldots,X_d]\!]$ which is a regular local 
ring of dimension $d=\dim(R)$. 
By the proof of~\cite{matsumura:crt} Theorem 29.4(iii), the map 
$R'\rightarrow R$ given by $X_i\mapsto x_i$ is injective and makes 
$R$ into a finite $R'$-module.  The fact that $x_1,\ldots,x_d$ 
generate a reduction of $\m$ implies, as in the proof of 
Theorem~\ref{thm:main}, that $\rank_{R'}(R)=e(R)=e(R_{\p})$ so that 
the proof is now identical to that of Theorem~\ref{thm:main}.
\end{proof}

In Example~\ref{ex:mixed} we show that the construction of Step 5 
does not work if the assumption ``$p$ is part of a reductive system 
of parameters of $R/\ann(M)$'' is dropped.

\section{The Conjecture of Kurano and Roberts} \label{sec:orig}

At this point, it seems wise to consider the status of 
Conjecture~\ref{conj:1} in light of Theorems~\ref{thm:main} 
and~\ref{thm:mixed}.  Conjecture~\ref{conj:1} holds for regular local 
rings containing a field by Theorem~\ref{thm:main}.  Of course, 
because the positivity conjecture holds for regular local rings 
containing a field, this also follows from Theorem~\ref{thm:0}.
In the mixed characteristic, unramified case, Theorem~\ref{thm:mixed} 
does not completely resolve this conjecture, but it shows us exactly 
where to focus our attention.  

Assume that $(R,\m)$ is an unramified regular local ring of mixed 
characteristic $p$ with prime ideals $\p$ and $\q$ such that 
$\sqrt{\p+\q}=\m$ and $\dim(R/\p)+\dim(R/\q)=\dim(R)$.  We want to 
show that $\p^{(n)}\cap\q\subseteq\m^{n+1}$.  Using standard methods, 
we may assume that $R$ is complete with infinite residue field.  
(Note that, for the question of Kurano and Roberts, we do not need to 
assume that $\p$ is analytically unramified to make this reduction.)  
Suppose that $f\in\p^{(n)}\cap\q$ and $f\not\in\m^{n+1}$.  Using the 
Associativity Formula, we may assume that $f$ is irreducible.  Let 
$R_1=R/fR$, $\m_1=\m R_1$, and so on.  The data $R_1,\m_1,\p_1,\q_1, 
M_1=R_1$ then gives a counterexample to Conjecture~\ref{conj:main} 
where the ring is in fact a complete domain.  By 
Theorems~\ref{thm:main} and~\ref{thm:mixed}, we therefore know that 
$R_1$ does not contain a field and the residual characteristic $p$ is 
not part of a reductive system of parameters of $R_1$.  An 
appropriate choice of variables then shows that we can write $R_1$ in 
the form
\[ R_1=V[\![X_1,\ldots,X_d]\!]/(p^n+a_1p^{n-1}+\cdots+a_n) \]
where $V$ is a complete $p$-ring, $a_i\in(X_1,\ldots,X_d)^i$ and 
$a_n\in\m^{n+1}$.  Therefore, with an appropriate choice of variables 
for $R$, we know that the only form $f$ can have is 
$f=p^n+a_1p^{n-1}+\cdots+a_n$.

\section{Examples} \label{sec:ex}

The following example demonstrates that, in 
Conjecture~\ref{conj:main}, the ring $R/\ann(M)$ must be
equidimensional.

\begin{ex} \label{ex:eq}
Let $R=M=k[\![X,Y,Z]\!]/(XY,XZ)=k[\![x,y,z]\!]$ and let $\p=(x)R$ and 
$\q=(y,z)R$.  Then
$R$ is excellent, $\p$ and $\q$ are in the support of $M$, $\p$ is 
analytically unramified,
$e(M_{\p})=1=e(M)$, and $\sqrt{\p+\q}=\m$.
However, $R$ is not equidimensional and
\[ \dim(R/\p)+\dim(R/\q)=3>2=\dim(M).\]
\end{ex}

The following example demonstrates that, in 
Conjecture~\ref{conj:main}, the condition 
$e(M_{\p})=e(M)$ is
necessary.

\begin{ex} \label{ex:mult}
Let $R=M=k[\![X,Y]\!]/(XY)=k[\![x,y]\!]$ and let $\p=(x)R$ and 
$\q=(y)R$.  Then
$R$ is excellent and equidimensional, $\p$ and $\q$ are in the 
support of $M$, $\p$ is analytically unramified,
and $\sqrt{\p+\q}=\m$.
However, $e(M_{\p})=1<2=e(M)$ and
\[ \dim(R/\p)+\dim(R/\q)=2>1=\dim(M).\]
\end{ex}

The following example demonstrates that, in 
Conjecture~\ref{conj:main}, the condition $\sqrt{\p+\q}=\m$ is
necessary.

\begin{ex} \label{ex:sqrt}
Let $R=M=k[\![X]\!]$ and let $\p=\q=(0)R$.  Then 
$R$ is excellent and equidimensional, $\p$ and $\q$ are in the 
support of $M$, $\p$ is analytically unramified,
and $e(M_{\p})=1=e(M)$.
However, $\sqrt{\p+\q}\neq\m$ and
\[ \dim(R/\p)+\dim(R/\q)=2>1=\dim(R). \]
\end{ex}

The following example demonstrates that, in 
Conjecture~\ref{conj:main}, the ideal $\q$ must be in 
the support of $M$.  (Note that the question
makes no sense if $\p$ is not in the support of $M$, since then 
$e(M_{\p})=0<e(M)$.)

\begin{ex} \label{ex:supp}
Let 
\[ R=k[\![ X,Y,Z,W]\!]/(XY,YZ,ZW,WX)=k[\![x,y,z,w]\!].\]  
Let $M=R/(x,z)\cong k[\![Y,W]\!]$.
Let $\p=(x,z)$ or $\p=(x,y,z)$ so that $e(M)=e(M_{\p})=1$ and 
$\dim(R/\p)\geq 1$.  Let
$\q=(y,w)$ so that $M_{\q}=0$ and $\dim(R/\q)=2$.  
Then $R$ is excellent, $R/\ann(M)$ is a domain, $\p\in\supp(M)$, $\p$ 
is analytically unramified,
and 
$\sqrt{\p+\q}=\m$.  However, $\q\not\in\supp(M)$ and
\[ \dim(R/\p)+\dim(R/\q)\geq 3>2=\dim(M). \]
\end{ex}

The following example demonstrates that, in 
Conjecture~\ref{conj:main}, the quotient ring $R/\ann(M)$ must 
be catenary.

\begin{ex} \label{ex:morenagata}
(\cite{nagata:lr} Appendix (E.2))
Let $K$ be a field and let $x$ be an indeterminate.  Let $r$ be a 
positive integer and
$z_1,\ldots,z_r\in K[\![x]\!]$ be power series which are 
algebraically independent over
$K(x)$.  Write $z_i=\sum_j a_{ij} x^j$ and for $j>0$ let
\[ z_{ij}=\frac{z_i-\sum_{k<j}a_{ik}x^k}{x^{j-1}}=\sum_{k\geq 
j}a_{ik}x^{k-j+1}. \]
Let $m$ be a positive integer and $y_1,\ldots,y_m$ be algebraically 
independent elements 
(indeterminates) over $K[x,z_1,\ldots,z_r]$.  Let 
$R_0=K[x,\{z_{ij}\}]$ and for $i=1,\ldots,m$
let $R_i=R_0[y_1,\ldots,y_i]$.  

The ideal $xR_0$ is maximal in $R_0$, because 
$z_{ij}=xz_{ij+1}+a_{ij}x\in xR_0$.  As Nagata
notes, $(R_0)_{xR_0}$ is a discrete valuation ring.  Let 
$\m_i=(x,y_1,\ldots,y_i)R_i$ and
$\n_i=(x-1,z_1,\ldots,z_r,y_1,\ldots,y_i)R_i$ for $i=0,\ldots,m$.  
The rings $V_i=(R_i)_{\m_i}$
are regular local rings of dimension $i+1$ and the rings 
$W_i=(R_i)_{\n_i}$ are regular local
rings of dimension $r+i+1$.  Let 
$S_i=R_i\smallsetminus(\m_i\cup\n_i)$ and let $A_i=(R_i)_{S_i}$.
The maximal ideals of $A_i$ are $\m_iA_i$ and $\n_iA_i$, and 
$(A_i)_{\m_iA_i}=V_i$ and
$(A_i)_{\n_iA_i}=W_i$.  It follows (\cite{nagata:lr} Appendix (E1.2)) 
that each $A_i$ is Noetherian
and $A_i/\m_i A_i=A_i/\n_i A_i=K$.

Let $A=A_m$ and let $J$ denote the Jacobson radical of $A$.  
Let $R=K+J$ which is a subring of $A$.  In fact,
by~\cite{nagata:lr} Appendix (E2.1) $R$ is a local ring with maximal 
ideal $J$, the residue field of $R$ is
$K$, and $A$ is a finite $R$-module.  We claim that $R$ is the 
desired example.  

First, we show that, for every prime ideal $\p$ of $R$, 
$e(R_{\p})=1$.  To show this, it suffices to 
show the following:  (i) $e(R)=1$ and (ii) for every nonmaximal prime 
ideal $\p$ of $R$, there
exists a unique prime ideal $P$ of $A$ such that $P\cap R=\p$ and 
that $R_{\p}=A_P$.  This
is sufficient because $A_P$ is a regular local ring.  
For (i) we note that Nagata's computation shows that
\[ e(R)=e((A_m)_{\n_mA_m})=e(W_m)=1. \]
For (ii), it suffices to verify the following.

Claim:  For every nonzero element $\lambda$ of the $R$-module $A/R$, 
$\ann_R(\lambda)=J$.

Before we prove the claim , we show how it implies (ii).  Let $\p$ be 
a non-maximal
prime ideal of $R$ and let $P$ be a prime ideal of $A$ contracting to 
$\p$ in $R$.  (Note that
the fact that the extension $R\rightarrow A$ is integral implies that 
such a prime $P$ always 
exists and is not maximal.)  
To verify (ii) it suffices to show that $(A/R)_{\p}=(A/R)\otimes_R 
R_{\p}=0$ because then we will 
have an isomorphism
$R_{\p}\cong A_{\p}$ so that $A_{\p}$ is a local ring;  the 
integrality of the extension 
$R\rightarrow A$
implies that there are no containments between prime ideals 
contracting to $\p$ in $R$ so there 
is a 
unique such prime in $A$.  Furthermore, it is straightforward to show 
that $A_P\cong A_{\p}$.  
By the claim, for every nonzero 
element $\lambda\in A/R$, $\ann_R(\lambda)=J$, and this implies that 
every such $\lambda$ is
annihilated by an element $s\in J\smallsetminus \p$.   If 
$\lambda_1,\ldots,\lambda_u\in A/R$
is a generating set of $A/R$ as an $R$-module, then 
$\lambda_1/1,\ldots,\lambda_u/1\in (A/R)_{\p}$
is a generating set of $(A/R)_{\p}$ over $R_{\p}$.  The fact that 
every $\lambda_i$ is annihilated by
an element $s_i\in J\smallsetminus \p$ implies that 
$\lambda_i/1=s_i\lambda_i/s_i=0$ so that
$(A/R)_{\p}=0$ as desired.

Now, we prove the claim.  If $t\in A\smallsetminus R$, let $\ol{t}$ 
denote the class of $t$ in the
quotient $A/R$.  The fact that $\ol{t}\neq 0$ implies that 
$\ann_R(\ol{t})\subseteq J$.  However,
$Jt\subseteq J\subset R$ which implies that $J\ol{t}=0$ so that 
$J\subseteq\ann(\ol{t})\subseteq J$
giving the desired equality.

As an ideal of $A$, $J$ is generated by the elements
$x(x-1),\{y_i\},\{z_j\}$.  As an ideal of $R$, 
\[ 
J=(x(x-1),\{y_i\},\{z_j\},\{xy_i\},\{xz_j\},\{(x-1)y_i\},\{(x-1)z_j\})R 
.\]  
Let $P=(x)A$ and $Q=(x-1)A$ which are prime ideals in $A$, and let 
$\p=P\cap R$ and $\q=Q\cap R$.
Then
\[ \p =(x(x-1),xy_1,\ldots,xy_m,xz_1,\ldots,xz_r)R \]
and
\[ \q =(x(x-1),(x-1)y_1,\ldots,(x-1)y_m,(x-1)z_1,\ldots,(x-1)z_r)R \]
so that $\p+\q=J$.  The condition $e(R_{\p})=e(R)$ is automatic.  
Furthermore, $\p$ is analytically unramified, as follows.  We have an 
inclusion of rings $R/\p\hookrightarrow A/P=A/xA$ such that the 
extension is module finite.  Also, $A/xA$ is a regular local ring and 
is therefore analytically unramified.  The maximal ideal of $R/\p$ 
extends to the maximal ideal of $A/xA$ so that the completion of 
$A/xA$ contains the completion of $R/\p$ (c.f., Atiyah and 
MacDonald~\cite{atiyah:ica} Theorem 10.11 and Corollary 10.3).  
Therefore, $\widehat{R/\p}$ is a domain.

Finally,
\[ \dim(R/\p)+\dim(R/\q)=\dim(A/P)+\dim(A/Q)=m+(r+m) \]
which is strictly greater than $\dim(R)=r+m+1$ if we assume that 
$m>1$.
\end{ex}

I do not know of an example where $R/\ann(M)$ is a local, catenary,
equidimensional ring which is not universally catenary and where 
Conjecture~\ref{conj:main}
does not hold.   

The following example demonstrates that, for a given ring $R$ and 
module $M$, the equality $e(M_{\p})=e(M)$ 
is weaker than the equality $e(R_{\p})=e(R)$, even when $R$ is 
complete, 
equidimensional and $\dim(M)=\dim(R)$.

\begin{ex} \label{ex:more}
Let 
\[ R=k[\![ X,Y,Z,W]\!]/(X,Y)\cap(Z,W)=k[\![x,y,z,w]\!].\]  
Let $M=R/(x,y)\cong k[\![Z,W]\!]$.  Then, $e(R)=2$ and, for every 
nonmaximal prime ideal $\p$ 
of $R$, $e(R_{\p})=1$.   However, for every prime $\p$ in the support 
of $M$, $e(M_{\p})=1=e(M)$.
\end{ex}

The following example demonstrates that, if $\p$ is a prime ideal in 
a local domain which is 
analytically unramified and does not satisfy the condition 
$\Ht(\p)+\dim(R/\p)=\dim(R)$, then 
the conclusion of Theorem~\ref{thm:loc} does not hold.

\begin{ex} \label{ex:nagata1}
We continue with the notation of Example~\ref{ex:morenagata}.  
Consider the polynomial
$f=y_m^s+x^t(x-1)\in A_{m-1}[y_m]$ where $s\geq t>1$.  
The fact that $A_{m-1}$ is a regular semilocal domain implies 
that $A_{m-1}$
is a unique factorization domain (\cite{nagata:lr} (28.8)).  The 
ideal $x$ generates a prime
ideal in $A_{m-1}$ and it follows from Eisenstein's criterion that 
$f$ is irreducible in
the polynomial ring $A_{m-1}[y_m]$.  Therefore, $f$ generates a prime 
ideal in $A_{m-1}[y_m]$.
Because $A_m$ is a localization of $A_{m-1}[y_m]$ and $f$ is a 
nonzero nonunit in $A_m$, we
see that the ideal $f A_m$ is prime.  Let $A=A_m/f A_m$ which is a 
semilocal domain with
maximal ideals $\m_m A$ and $\n_m A$.  

Let $J$ denote the Jacobson radical of $A$ and let $R=K+J$ which is a 
subring of $A$.  As
in Example~\ref{ex:morenagata} $R$ is a local ring with maximal ideal 
$J$, the residue field 
of $R$ is
$K$, and $A$ is a finite $R$-module.  We claim that $R$ is the 
desired example.  

First, we verify that $e(R)=1$.  This is straightforward, as Nagata's 
computation shows that
\[ e(R)=e((A_m)_{\n_mA_m}/(f))=e(W_m/(f)). \]
The fact that $x-1$ is a minimal generator of the maximal ideal of 
$W_m$ and that $x$ is a unit
in $W_m$ implies that $f=y_m^s+x^t(x-1)$ is a minimal generator of 
the maximal ideal of $W_m$.
Because $W_m$ is a regular local ring, it follows that $e(R)=1$.

In order to verify that there is a prime ideal $\p$ of $R$ such that 
$e(R_{\p})>1$, it suffices to verify the
following.
Claim:  If $\p$ is a non-maximal prime ideal of $R$, then there 
exists a unique prime 
ideal $P$ of $A$ contracting to $\p$ in $R$ and $A_P=R_{\p}$.  The 
verification of this claim
is identical to that in Example~\ref{ex:morenagata}, so we omit it 
here.

We show how the claim implies the desired result.  Let 
$P=(x,y_m)A$.  Then $P$ is a prime ideal of $A$ because $x$ is a 
prime element of $A_{m-1}$ and
therefore the ideal $(x,y_m)A_m$ is a prime ideal of $A_m$ containing 
$f$.
Let $\p=P\cap R$.  Then $\p$ is not maximal because the extension 
$R\rightarrow A$ is integral 
and $P$ is not maximal.  The verification of the fact that $\p$ is 
analytically unramified is 
similar to that in
Example~\ref{ex:morenagata}.  (By Theorem~\ref{thm:loc} it follows 
that $\Ht(\p)+\dim(R/\p)\neq\dim(R)$.)
As noted, $R_{\p}\cong A_P$ and it follows that
\[ e(R_{\p})=e(A_P)=e(V_m/fV_m). \]
The fact that $x$ is a minimal generator of the maximal ideal of 
$V_m$ and that $x-1$ is a unit
in $V_m$ implies that $f\in (\m_m V_m)^t\smallsetminus(\m_m 
V_m)^{t+1}$.  The fact that $V_m$
is a regular local ring implies that $e(R_{\p})=t>1$, as desired.
\end{ex} 

I do not know of a local, catenary,
equidimensional ring which is not universally catenary where this 
inequality of multiplicities
does not hold.  A number of examples have been constructed (e.g., 
\cite{nagata:lr} Appendix (E.2), 
Heinzer, Rotthaus and Wiegand~\cite{heinzer:ndihic} Example (4.5), 
Nishimura~\cite{nishimura:felrI} 
Sections 2,4 and 6) where the hypotheses of Theorem~\ref{thm:loc} do 
not hold.  However, in each 
example, the conclusion of Theorem~\ref{thm:loc} does hold.

The following example demonstrates that if $R$ is a complete domain of
mixed characteristic $p$ such that $p$ is not
part of a reductive system of parameters of $R$, then there does not
exist in general a regular local ring $R'\subseteq R/\ann(M)$ such 
that
the extension is finite and $e(R)=\rank_R'(R)$.  In particular, in 
Theorem~\ref{thm:mixed}, if the assumption ``$p$ is part of a 
reductive system of parameters of $R/\ann(M)$'' is dropped, then the 
construction of Step 5 of the proof does not work.
 
\begin{ex} \label{ex:mixed}
Let $p$ be a prime number and $(V,pV,k)$ a complete $p$-ring.  Let 
$A=V[\![X,Y]\!]$, $f=pX+Y^3$ and 
$R=A/(f)$.  It is straightforward to verify that $f$ is irreducible 
in $A$ so that $R$ is a 
complete domain of dimension 2 and mixed characteristic $p$.  Let 
$\m$ and $\n$ denote the maximal 
ideals of $R$ and $A$, respectively.

The element $p$ is not part of a reductive system of parameters of 
$R$, as follows.  Suppose the 
contrary.  Because $\dim(R)=2$ there would exist an element $z\in R$ 
such that $p,z$ is a reductive system of parameters of $R$.  Let 
$\q=(p,z)R$.  
By~\cite{northcott:rilr} Theorem 4.2, $\q$ is a minimal reduction of 
$\m$.
By Herrmann, Ikeda and Orbanz~\cite{herrmann:ebu} Proposition 10.17, 
the fact that $p,z$ generate a minimal reduction of $\m$ implies that 
the initial forms of $p$ and $z$ in the associated graded ring 
$\mathrm{gr}_{\m}(R)$ form a homogeneous system of parameters for 
$\mathrm{gr}_{\m}(R)$.  The fact that $A$ is regular and $R=A/(f)$ 
implies that $\mathrm{gr}_{\m}(R)$ is the quotient of 
$\mathrm{gr}_{\n}(A)$ by the initial form of $f$.  More specifically, 
let $\mathbf{P}, \mathbf{X}, \mathbf{Y}$ denote the initial forms of 
$p,X,Y$ in $\n/\n^2$, respectively.  Then $\mathbf{P}, \mathbf{X}, 
\mathbf{Y}$ are indeterminates in $\mathrm{gr}_{\n}(A)$ and 
$\mathrm{gr}_{\n}(A)$ is the polynomial ring 
$k[\mathbf{P},\mathbf{X},\mathbf{Y}]$.  Then $\mathrm{gr}_{\m}(R)$ is 
a quotient of $\mathrm{gr}_{\n}(A)$:
\[ 
\mathrm{gr}_{\m}(R)=k[\mathbf{P},\mathbf{X},\mathbf{Y}]/(\mathbf{PX}). 
\]
The initial form of $p$ in $\mathrm{gr}_{\m}(R)$ is the image of 
$\mathbf{P}$ which is contained in a minimal prime ideal of 
$\mathrm{gr}_{\m}(R)$.  Therefore, this initial form can not be part 
of a homogeneous system of parameters of $\mathrm{gr}_{\m}(R)$, 
giving the desired contradiction.

Now, we show that there does not exist a complete regular local ring 
$R'$ contained in $R$ such that the extension $R'\rightarrow R$ is 
module-finite and $e(R)=\rank_{R'}(R)$.  This will prove that we have 
the desired example.  Suppose that such a ring $R'$ existed.  Because 
$R$ is a domain of mixed characteristic $p$, $R$ does not contain a 
field.  Therefore $R'$ does not contain a field.  Because $R$ is 
unramified, $R'$ is also unramified because otherwise, 
$p\in(\m')^2\subseteq\m^2$ which contradicts the fact that $R$ is 
unramified.  Therefore, $R'$ is of the form $W[\![Z]\!]$ where 
$(W,pW)$ is a complete $p$-ring and $Z$ is analytically independent 
over $W$.  By assumption, $e(\m,R)=e(R)=\rank_{R'}(R)=e(\m'R,R)$.  By 
the theorem of Rees (c.f., \cite{herrmann:ebu} Theorem 19.5) this 
implies that $\m' R$ is a reduction of $\m$.  Because $\m' R=(p,Z)R$ 
we see that this implies that $p$ is part of a reductive system of 
parameters of $R$, a contradiction.

We also note that it is straightforward to verify that no change a 
variables will remedy the noted behavior of this example.  
\end{ex}

\bigskip
\Large
\noindent \textbf{Acknowledgements}  
\normalsize

\medskip
I am grateful to Kazuhiko Kurano for suggesting the present version 
of Lemma~\ref{lem:1}
which is more general than the previous version.
I am also grateful to Dan Katz for many stimulating discussions 
regarding these 
conjectures and theorems.   Finally, I would like to thank the 
Department of 
Mathematics  at the University of Kansas at Lawrence for its 
hospitality during the preparation of this paper.


\providecommand{\bysame}{\leavevmode\hbox to3em{\hrulefill}\thinspace}

\end{document}